# The Stein hull

**Clément Marteau**

*Université de Provence, CMI,*
*39, rue F. Joliot-Curie,*
*F-13 453 Marseille cedex 13, France.*
**e-mail:** clement.marteau@cmi.univ-mrs.fr

**Abstract:** We are interested in the statistical linear inverse problem $Y = Af + \epsilon\xi$, where $A$ denotes a compact operator and $\epsilon\xi$ a stochastic noise. In a first time, we investigate the link between some threshold estimators and the risk hull point of view introduced in (5). The penalized blockwise Stein's rule plays a central role in this study. In particular, this estimator may be considered as a risk hull minimization method, provided the penalty is well-chosen. Using this perspective, we study the properties of the threshold and propose an admissible range for the penalty leading to accurate results. We eventually propose a penalty close to the lower bound of this range. The risk hull point of view provides interesting tools for the construction of adaptive estimators. It sheds light on the processes governing the behavior of linear estimators. The variability of the problem may be indeed quite large and should be carefully controlled.

## Contents

## 1. Introduction

This paper deals with the statistical inverse problem:

$$Y = Af + \epsilon\xi, \tag{1.1}$$

where $H, K$ are Hilbert spaces and $A : H \to K$ denotes a linear operator. The function $f \in H$ is unknown and has to be recovered from a measure of $Af$ corrupted by some stochastic noise $\epsilon\xi$. Here, $\epsilon$ represents a positive noise level and $\xi$ a Gaussian white noise (see (14) for more details). In particular, for all $g \in K$, we can observe:

$$\langle Y, g \rangle = \langle Af, g \rangle + \epsilon \langle \xi, g \rangle, \tag{1.2}$$

where $\langle \xi, g \rangle \sim \mathcal{N}(0, \|g\|^2)$. Denote by $A^\star$ the adjoint operator of $A$. In the sequel, $A$ is supposed to be a compact operator. Such a restriction is rather interesting from a mathematical point of view. The operator $(A^\star A)^{-1}$ is unbounded: the least square solution $\hat{f}_{LS} = (A^\star A)^{-1} A^\star Y$ does not continuously depend on $Y$. The problem is said to be ill-posed.

Several studies of ill-posed inverse problems in a statistical context were proposed in recent years. It would be however impossible to cite them all. For the interested reader, we may mention (12) and (11) for convolution operators, (15) for the positron emission tomography problem, (9) in a wavelett setting, or (2) for a general statistical approach and some rates of convergence. We refer also to (10) for a survey in a numerical setting.

Using a specific representation (i.e. particular choices for $g$ in (1.2)) may help the understanding of the model (1.1). In this sense, the classical singular value decomposition (SVD) is a very useful tool. Since $A^\star A$ is compact and auto-adjoint, the associated sequence of eigenvalues $(b_k^2)_{k\in\mathbb{N}}$ is strictly positive and converges to 0 as $k \to +\infty$. The sequence of eigenvectors $(\phi_k)_{k\in\mathbb{N}}$ is supposed in the sequel to be an orthonormal basis of $H$. For all $k \in \mathbb{N}$, set $\psi_k = b_k^{-1} A\phi_k$. The triple $(b_k, \phi_k, \psi_k)_{k\in\mathbb{N}}$ verifies:

$$\begin{cases} A\phi_k = b_k \psi_k, \\ A^\star \psi_k = b_k \phi_k, \end{cases} \tag{1.3}$$

for all $k \in \mathbb{N}$. This is the singular value decomposition of $A^\star A$. This representation leads to a simpler understanding of the model (1.1). Indeed, for all $k \in \mathbb{N}$, using (1.3) and the properties of the Gaussian white noise:

$$y_k = \langle Y, \psi_k \rangle = \langle Af, \psi_k \rangle + \epsilon \langle \xi, \psi_k \rangle = b_k \langle f, \phi_k \rangle + \epsilon \xi_k, \tag{1.4}$$

where the $\xi_k$ are i.i.d. standard Gaussian variables. Hence, for all $k \in \mathbb{N}$, we can obtain from (1.1) an observation on $\theta_k = \langle f, \phi_k \rangle$. In the $\ell^2$-sense, $\theta = (\theta_k)_{k\in\mathbb{N}}$ and $f$ represent the same mathematical object. The sequence space model (1.4) clarifies the effect of $A$ on the signal $f$. Since $A$ is compact, $b_k \to 0$ as $k \to +\infty$. For large values of $k$, the coefficients $b_k\theta_k$ are negligible compared to $\epsilon\xi_k$. In a certain sense, the signal is smoothed by the operator. The recovering becomes difficult in the presence of noise for large 'frequencies', i.e. when $k$ is large.

From now, our aim is to estimate the sequence $(\theta_k)_{k\in\mathbb{N}}$. The linear estimation plays an important role in the inverse problem framework and is a starting point for several recovering methods. Let $(\lambda_k)_{k\in\mathbb{N}}$ be a real sequence with values in $[0, 1]$. In the following, this sequence will be called a filter. The associated linear estimator is defined by:

$$\hat{f}_\lambda = \sum_{k=1}^{+\infty} \lambda_k b_k^{-1} y_k \phi_k.$$

In the sequel, $\hat{f}_\lambda$ may be sometimes identified with $\hat{\theta}_\lambda = (\lambda_k b_k^{-1} y_k)_{k\in\mathbb{N}}$. The meaning will be clear following the context. The error related to $\hat{f}_\lambda$ is measured through the quadratic risk:

$$R(\theta,\lambda) = \mathbb{E}_\theta \|\hat{f}_\lambda - f\|^2 = \sum_{k=1}^{+\infty}(1-\lambda_k)^2\theta_k^2 + \epsilon^2 \sum_{k=1}^{+\infty} \lambda_k^2 b_k^{-2} = \mathbb{E}_\theta\|\hat{\theta}_\lambda - \theta\|^2. \quad (1.5)$$

The most interesting filters are data-driven in the sense that they do not require a priori informations on $f$. In order to evaluate the performances of a data-driven estimator, we use in this paper the oracle point of view. Given a family of estimators $T$, define:

$$\theta_T = \arg\inf_{\tilde{\theta}\in T} \mathbb{E}_\theta\|\tilde{\theta} - \theta\|^2.$$

This is the oracle for the family T: it is the best possible estimator of $\theta$ in the family $T$. A data-driven estimator $\theta^\star$ may be compared to $\theta_T$ via the oracle inequality:

$$\mathbb{E}_\theta\|\theta^\star - \theta\|^2 \leq (1+\vartheta_\epsilon)\mathbb{E}_\theta\|\theta_T - \theta\|^2 + C\epsilon^2, \quad (1.6)$$

with $\vartheta_\epsilon, C > 0$. The quantity $C\epsilon^2$ is a residual term. The inequality (1.6) is said to be sharp if $\vartheta_\epsilon \to 0$ as $\epsilon \to 0$: $\theta^\star$ asymptotically mimics the behaviour of $\theta_T$. Oracle inequalities play an important, though recent role in statistics. They provide a precise measure on the performances of $\theta^\star$. They do not require a priori informations on the signal and are non-asymptotic. In several situations, oracle results may also lead to interesting minimax rates of convergence. This theory has given rise to a considerable amount of literature. We mention in particular (9), (1), (6) or (3) for a survey.

There exist several approaches leading to accurate oracle inequalities. The unbiased risk estimation (URE) method presents an interesting behavior. However, it does not take into account the variability of the problem. This is quite problematic in the inverse problem framework. The risk hull minimization (RHM) method initiated in (5) is an interesting alternative. This method proposes a data-driven bandwidth for projection schemes in the SVD setting. The principle is to identify the stochastic processes that control the behavior of a projection estimator. Then, a deterministic criterion, called a hull, is constructed in order to contain these processes. In several cases, this approach leads to an accurate recovering.

Our aim in this paper is to obtain oracle inequalities on wide families of estimators. In order to achieve this goal, we will consider the penalized blockwise Stein's rule initiated by (8). Some specific choices of penalty have already been proposed. Here, we present a general approach. The risk hull point of view may precise the role of the penalty. In particular, we present a link between the RHM procedure and such thresholds estimators.

This paper is organized as follows. The construction of the penalized blockwise Stein's rule and some related results are recalled in Section 2. Section 3 presents a link between hulls and such threshold estimators. Section 4 investigates the performances of the penalized blockwise Stein's estimator following the chosen penalty. Section 5 proposes some examples and a discussion on the choice of the penalty. Some results on the theory of ordered processes and the proofs of the main results are eventually gathered in Section 6.

## 2. The penalized blockwise Stein's rule

The construction of adaptive estimators is an interesting problematic. In the oracle sense, an ideal goal of adaptation is to obtain a sharp oracle inequality over all possible estimators. This is in most cases an unreachable task since this set is rather large. The difficulty of the oracle adaptation increases with the size of the considered family. In this paper, we restrict ourselves to the linear and monotone estimators. In the following, this family will be identified to the collection:

$$\Lambda_{mon} = \left\{ \lambda = (\lambda_k)_{k\in\mathbb{N}} \in \ell^2 : 1 \geq \lambda_1 \geq \cdots \geq \lambda_k \geq \cdots \geq 0 \right\},$$

of linear and monotone filters. This family contains most of the existing linear procedures. We may mention for instance the spectral cut-off, Tikhonov, Pinsker or the Landweber filters (see for instance (10) or (2)).

A good way to obtain oracle inequalities on $\Lambda_{mon}$ is to consider in a first time the family of blockwise constant filters $\Lambda^\star$ defined by:

$$\Lambda^\star = \left\{ \lambda \in l^2 : \ 0 \leq \lambda_k \leq 1, \ \lambda_k = \lambda_{K_j}, \ \forall k \in [K_j, K_{j+1}-1], \right. \\ \left. j = 0, \dots, J, \lambda_k = 0, \ k > N \right\}, \tag{2.1}$$

where $J$, $N$ and $(K_j)_{j=0\dots J}$ are such that $K_0 = 1$, $K_J = N+1$ and $K_j > K_{j-1}$. In the following, we will also use the notations $I_j = \{k \in [K_{j-1}, K_j - 1]\}$ and $T_j = K_j - K_{j-1}$, for all $j \in \{1, \dots, J\}$.

The family $\Lambda^\star$ can easily be handled. In particular, each block $I_j$ can be studied independently of the other ones. This simplifies considerably the study of the considered estimators. Moreover, for all $\theta \in \ell^2$,

$$R(\theta, \lambda_{mon}) = \inf_{\lambda \in \Lambda_{mon}} R(\theta, \lambda) \text{ and } R(\theta, \lambda^0) = \inf_{\lambda \in \Lambda^\star} R(\theta, \lambda), \tag{2.2}$$

are in fact rather close, subject to some reasonnable constraints on the sequences $(b_k)_{k\in\mathbb{N}}$ and $(T_j)_{j=1\dots J}$ (see Section 5 or (8) for more details).

In order to construct a data-driven filter on $\Lambda^\star$, one may consider the well-known unbiased risk estimation approach (URE). The principle is rather intuitive. We want to construct a filter as close as possible to the oracle on $\Lambda^\star$. The quadratic risk $R(\theta, \lambda)$ associated to each filter is unknown: it explicitly depends on the function $f$, i.e. the sequence $\theta$. This term can however be approximated by an estimator $U(y,\lambda)$. The related adaptive filter is then defined as:

$$\begin{aligned}\hat{\lambda}^{URE} &= \arg\inf_{\lambda\in\Lambda^\star} U(y,\lambda),\\ &= \arg\inf_{\lambda\in\Lambda^\star}\left[\sum_{k=1}^{+\infty}(\lambda_k^2-2\lambda_k)(b_k^{-2}y_k^2-\epsilon^2 b_k^{-2})+\epsilon^2\sum_{k=1}^{+\infty}\lambda_k^2 b_k^{-2}\right], \quad (2.3)\end{aligned}$$

and can explicitly be computed:

$$\hat{\lambda}_k^{URE} = \begin{cases}\left(1-\frac{\sigma_j^2}{\|\tilde{y}\|_{(j)}^2}\right)_+ & , \quad k\in I_j,\ j=1\dots J,\\ 0 & , \quad k>N,\end{cases}$$

with,

$$\|\tilde{y}\|_{(j)}^2 = \epsilon^2\sum_{k\in I_j} b_k^{-2}y_k^2, \text{ and } \sigma_j^2=\epsilon^2\sum_{k\in I_j} b_k^{-2},\ \forall j\in\{1\dots J\}. \quad (2.4)$$

For all $j\in\{1,\dots,J\}$, $\sigma_j^2$ is in fact compared to $\|\tilde{y}\|_{(j)}^2$. Recall that $\mathbb{E}_\theta\|\tilde{y}\|_{(j)}^2 = \sum_{k\in I_j}\theta_k^2+\sigma_j^2$. In a certain sense, $\hat{\lambda}^{URE}$ only keeps the blocks where the signal is not negligible compared to $\sigma_j^2$. However, $\hat{\lambda}^{URE}$ is rather unstable, due to the ill-posedness of the problem. It is only concerned with the average behavior of the loss, or equivalently of $\|\tilde{y}\|_{(j)}^2$. The URE approach does not take into account the variability of the problem: the variance of $\|\tilde{y}\|_{(j)}^2$ explodes for large $j$. A complete discussion illustrated by some numerical simulations is provided in (5) in a slightly different setting.

As a benchmark for the URE method limitations, one may consider penalized estimators. Let $(\text{pen}_j)_{j=1\dots J}$ a real positive sequence. For all $\lambda\in\Lambda^\star$, consider the penalized unbiased risk estimator of the quadratic risk:

$$\begin{aligned}U_p(y,\lambda) &= \sum_{j=1}^{J}\left[\{\lambda_{K_j}^2-2\lambda_{K_j}\}(\|\tilde{y}\|_{(j)}^2-\sigma_j^2)+\lambda_{K_j}^2\sigma_j^2+2\lambda_{K_j}\text{pen}_j\right],\\ &= U(y,\lambda)+2\sum_{j=1}^{J}\lambda_{K_j}\text{pen}_j, \quad (2.5)\end{aligned}$$

where $U(y,\lambda)$ is defined in (2.3). The penalty should ideally contain the variability of the problem. An adaptive estimator may be constructed as the minimizer of $U_p(y,\lambda)$ on the family $\Lambda^\star$. The solution is:

$$\lambda_k^\star = \begin{cases}\left(1-\frac{\sigma_j^2+\text{pen}_j}{\|\tilde{y}\|_{(j)}^2}\right)_+ & , \quad k\in I_j,\ j=1\dots J,\\ 0 & , \quad k>N.\end{cases} \quad (2.6)$$

Cavalier and Tsybakov (2002) proposed to choose $\text{pen}_j=\varphi_j\sigma_j^2$ for all $j\in\{1,\dots,J\}$. Concerning the sequence $(\varphi_j)_{j=1\dots J}$, they set the following condition:

**Assumption A1**: *There exists a constant $C_1$ independent of $\epsilon$ such that:*

$$\sum_{j=1}^{J} \max_{k\in I_j} b_k^{-2} \exp\left[ -\frac{\varphi_j^2}{16\Delta_j(1+2\sqrt{\varphi_j})}\right] \leq C_1, \ \textit{with}\ \Delta_j = \frac{\max_{k\in I_j} \epsilon^2 b_k^{-2}}{\sigma_j^2}, \quad (2.7)$$

*for all $j \in \{1,\dots,J\}$.*

The performances of the related estimator $\tilde{\theta}$ are summarized in the following theorem:

**Theorem 1.** *(Cavalier and Tsybakov (2002)).*
*Let $\tilde{\theta}$ the estimator associated to the filter $\lambda^\star$ defined in (2.6) with the penalty $\mathrm{pen}_j = \varphi_j \sigma_j^2$ for all $j \in \mathbb{N}$. Assume that the sequence $(\varphi_j)_{j\in 1\dots J}$ satisfies Assumption A1 and $\varphi_j \leq 1 - 4\Delta_j$, for all $j \in \{1,\dots,J\}$. Then, for all $\theta \in \ell^2$ and $0 < \epsilon < 1$, we have:*

$$\mathbb{E}_\theta \|\theta^\star - \theta\|^2 \leq (1+\varphi_\epsilon) \inf_{\lambda\in\Lambda^\star} R(\theta,\lambda) + 8C_1\epsilon^2,$$

*with $\varphi_\epsilon = \max_{1\leq j\leq J}(2\varphi_j + 16\Delta_j/\varphi_j)$.*

For a large range of inverse problems, Theorem 1 provides a sharp oracle inequality on $\Lambda^\star$ for $\lambda^\star$ with $\mathrm{pen}_j = \varphi_j\sigma_j^2$ for all $j \in \{1\dots J\}$. Indeed, assume for instance that the sequence $(b_k)_{k\in\mathbb{N}}$ possesses a polynomial decay. In this case, with an appropriate choice of blocks, $\Delta_j \to 0$ as $\epsilon \to 0$. It remains to choose $\varphi_j$ such that $\varphi_\epsilon \to 0$ as $\epsilon \to 0$.

In this paper, we investigate the relationship between the properties of the penalty $(\mathrm{pen}_j)_{j=1\dots J}$ and the performances of the related estimator. In particular, we shed some light on the link between such threshold estimators and the risk hull point of view introduced in (5). The underneath aim of this study is to precise the role of the penalty and provide an admissible range in the oracle sense.

## 3. Risk hull and penalties

The principle of risk hull minimization method has been introduced in (5). It provides an adaptive bandwidth choice for the projection (also called spectral cut-off) filters: $(\mathbf{1}_{\{k\leq N\}})_{k\in\mathbb{N}}$ with $N \in \mathbb{N}$. The projection estimation in the SVD formalism (1.4) may be seen as a toy model for adaptive estimation. Nevertheless, it presents several difficulties and requires a careful treatment. In this section, we recall the principle of the risk hull minimization for projection schemes. Then, we prove that the penalized quadratic risk is a hull for the family of blockwise constant filters. This result requires some conditions on the penalty $(\mathrm{pen}_j)_{j=1\dots J}$: the role of this sequence is then precised.

For all $N \in \mathbb{N}$, denote by $\hat{\theta}_N$ the projection estimator associated to the filter $(\mathbf{1}_{\{k \leq N\}})_{k\in\mathbb{N}}$. For each value of $N \in \mathbb{N}$, the related quadratic risk is:

$$\mathbb{E}_\theta \|\hat{\theta}_N - \theta\|^2 = \mathbb{E}_\theta \sum_{k=1}^{N} (b_k^{-1} y_k - \theta_k)^2 + \sum_{k>N} \theta_k^2 = \sum_{k>N} \theta_k^2 + \epsilon^2 \sum_{k=1}^{N} b_k^{-2}. \quad (3.1)$$

The optimal choice for $N$ is the oracle $N^0$ that minimizes $\mathbb{E}_\theta \|\hat{\theta}_N - \theta\|^2$. It is a trade-off between the two sums (bias and variance) in the r.h.s. of (3.1).

In order to construct an adaptive bandwidth $N^\star$, one may use the classical URE procedure. This approach has been studied in this setting by (6). The theoretical results are interesting but the numerical simulations may be somewhat disapointing in several cases. The URE approach does not take into account the variability of the problem, which may be quite large when considering compact operators.

The construction of an adaptive filter can be decomposed in two steps. First: evaluate the error associated to each filter by constructing an appropriate criterion. Then, use this criterion in order to propose an adaptive estimator. The criterion associated to the URE approach is the quadratic risk. It corresponds to the average behavior of the considered estimators and does not detect the variability. (5) are interested instead in the loss:

$$\|\hat{\theta}_N - \theta\|^2 = \sum_{k>N} \theta_k^2 + \epsilon^2 \sum_{k=1}^{N} b_k^{-2} \xi_k^2.$$

As a criterion, they use a deterministic term $V(\theta, N)$, called a hull, satisfying:

$$\mathbb{E}_\theta \sup_{N\in\mathbb{N}} [l(\theta, N) - V(\theta, N)] \leq 0. \quad (3.2)$$

This hull bounds uniformly the loss in the sense of inequality (3.2). It is constructed in order to contain the variability of the projection estimators.

The risk hull point of view makes the role of the stochastic processes involved in linear estimation more precise. In the same way, it quantifies the limitations of the URE method and may lead to an accurate understanding of the problem. The hull proposed by (5) precisely take into account the variability of the problem. The related adaptive estimator possesses interesting theoretical and numerical properties. In the same spirit, we mention (17) for general families of linear estimators.

The risk hull point of view may provide an interesting perspective on the blockwise constant adaptive approach. In this section, we propose sufficient conditions on the penalty making the penalized quadratic risk a hull.

First, we introduce some notations. For all $j \in \{1, \dots, J\}$, let $\eta_j$ defined by:

$$\eta_j = \epsilon^2 \sum_{k\in I_j} b_k^{-2} (\xi_k^2 - 1). \quad (3.3)$$

The random variable $\eta_j$ plays a central role in blockwise constant estimation. It corresponds to the main stochastic part of the loss in each block $I_j$. The hull proposed in Theorem 2 bellow is constructed in order to contain these terms. Introduce also:

$$\rho_\epsilon = \max_{j=1\ldots J} \sqrt{\Delta_j} \text{ and } \|\theta\|^2_{(j)} = \sum_{k\in I_j} \theta_k^2,\ \forall j \in \{1,\ldots,J\}, \tag{3.4}$$

where $\Delta_j$ is defined in (2.7). We will see that $\rho_\epsilon \to 0$ as $\epsilon \to 0$ with appropriate choices of blocks and minor assumptions on the sequence $(b_k)_{k\in\mathbb{N}}$ (see Section 5 for more details). Concerning the penalty, we set the following condition:

**Assumption A2**: *There exists a constant $C_2$ independent of $\epsilon$ such that:*

$$\sum_{j=1}^{J} \mathbb{E}\left[\eta_j - \mathrm{pen}_j\right]_+ \leq C_2\epsilon^2, \tag{3.5}$$

Assumptions A1 and A2 are in fact rather close. Indeed, the exponential term in A1 corresponds to the probability of the event $\{\lambda^\star_{K_j} > 0\}$ when the signal to noise ratio is small on the block $j$. In such a situation, $\|\tilde{y}\|^2_{(j)}$ is close to $\sigma_j^2 + \eta_j$. Recall that for all $j \in \{1\ldots J\}$,

$$\lambda^\star_{K_j} = \left(1 - \frac{\sigma_j^2 + \mathrm{pen}_j}{\|\tilde{y}\|^2_{(j)}}\right)_+ .$$

If the penalty is 'well-chosen', $P(\tilde{\lambda}_j > 0)$ is small, i.e. the penalty controls in a certain sense the variables $\eta_j$. This is exactly the principle of Assumption A2. In Section 4, the link between these two hypotheses will be strengthened via an upper bound of the l.h.s. of (3.5).

The proof of the following result is presented in Section 6.

**Theorem 2.** *Assume that Assumption A2 holds. Then, there exists $B > 0$ such that:*

$$\begin{aligned} V(\theta,\lambda) \;=\; (1+B\rho_\epsilon)\left\{\sum_{j=1}^{J}\left[(1-\lambda_{K_j})^2\|\theta\|^2_{(j)} + \lambda^2_{K_j}\sigma_j^2 + 2\lambda_{K_j}\mathrm{pen}_j\right] + \sum_{k>N}\theta_k^2\right\} \\ + C_2\epsilon^2 + B\rho_\epsilon R(\theta,\lambda^0), \end{aligned} \tag{3.6}$$

*is a risk hull on $\Lambda^\star$, i.e.:*

$$\mathbb{E}_\theta \sup_{\lambda\in\Lambda^\star}\left\{\|\hat{\theta}_\lambda - \theta\|^2 - V(\theta,\lambda)\right\} \leq 0.$$

Theorem 2 states in fact that the penalized quadratic risk:

$$R_p(\theta,\lambda) = \sum_{j=1}^{J}\left[(1-\lambda_{K_j})^2\|\theta\|^2_{(j)} + \lambda^2_{K_j}\sigma_j^2 + 2\lambda_{K_j}\mathrm{pen}_j\right] + \sum_{k>N}\theta_k^2, \tag{3.7}$$

is, up to some constants and residual terms, a risk hull on the family $\Lambda^\star$. Hence, $R_p(\theta,\lambda)$ may be a good criterion for the choice of $\lambda^\star$, provided that inequality (3.5) is satisfied. In a certain sense, Theorem 2 justifies the approach presented in Section 2 since the term $U_p(y,\lambda)$ defined in (2.5) is an estimator of $R_p(\theta,\lambda)$. The Theorem 2 makes the role of the penalty more precise. In order to obtain a hull, we need to construct a penalty that contains, in the sense of inequality (3.5), the variables $(\eta_j)_{j=1\dots J}$.

This result is established for general blocks. Following (7), there exist several choices that may lead to interesting results. An example is presented in Section 5.

The penalty $(\text{pen}_j)_{j=1\dots J}$ in the hull $V(\theta,\lambda)$ is associated in each block $I_j$ to the term $2\lambda_j$. The construction is then closely related to the problem of proving that:

$$\sum_{j=1}^J \mathbb{E}_\theta \sup_{\lambda_j\in[0,1]} \left\{\lambda_j^2\eta_j - 2\lambda_j \text{pen}_j\right\} \leq C_2\epsilon^2,$$

see the proof of Theorem 2 in Section 6. The penalty may also be associated to the term $\lambda_j^2$. In this case, we obtain the following result.

**Theorem 3.** *Assume that Assumption A2 holds. Then, there exists $B>0$ such that:*

$$W(\theta,\lambda) = (1+B\rho_\epsilon)\left\{\sum_{j=1}^J\left[(1-\lambda_{K_j})^2\|\theta\|^2_{(j)} + \lambda_{K_j}^2\sigma_j^2 + \lambda_{K_j}^2\text{pen}_j\right] + \sum_{k>N}\theta_k^2\right\} + C_2\epsilon^2 + B\rho_\epsilon R(\theta,\lambda^0), \quad (3.8)$$

*is a risk hull on $\Lambda^\star$, i.e.*

$$\mathbb{E}_\theta \sup_{\lambda\in\Lambda^\star}\left\{\|\hat\theta_\lambda-\theta\|^2 - W(\theta,\lambda)\right\} \leq 0.$$

Both hulls $V(\theta,\lambda)$ and $W(\theta,\lambda)$ contain the loss in the same way, i.e. the residual terms and the constant $B$ are exactly the same. However, $W(\theta,\lambda)$ is slighty smaller than $V(\theta,\lambda)$ since $0<\lambda_j<1$ for all $j\in\{1\dots J\}$. In some sense, $W(\theta,\lambda)$ is more precise. Nevertheless, we will use $V(\theta,\lambda)$ as a criterion in the following: the associated estimator has a simpler form. We will also see in Section 4 that the variability of the estimator of $V(\theta,\lambda)$ is contained by the penalty. This is not the case for $W(\theta,\lambda)$.

The proof of Theorem 3 is presented in Section 6. It follows essentially the same lines of the proof of Theorem 2.

## 4. Oracle inequalities

In Section 3, we have proposed a family of hull indexed by the penalty $(\text{pen}_j)_{j=1\dots J}$. A hull may be a good criterion in order to evaluate the performances of the estimators contained in $\Lambda^\star$ since it takes into account the variability of the problem. In this section, we are interested in the performances of the estimators constructed from these hulls. We are looking for conditions on the penalty $(\mathrm{pen}_j)_{j=1\dots J}$ that may lead to sharp oracle inequalities.

In the sequel, we set $\lambda_j = \lambda_{K_j}$ for all $j \in \{1\dots J\}$. This is a slight abuse of notation but the meaning will be clear following the context. Then define:

$$U_p(y,\lambda) = \sum_{j=1}^{J} \left[ (\lambda_j^2 - 2\lambda_j)(\|\tilde{y}\|_{(j)}^2 - \sigma_j^2) + \lambda_j^2 \sigma_j^2 + 2\lambda_j \mathrm{pen}_j \right].$$

The term $U_p(y,\lambda)$ is an estimator of the penalized quadratic risk $R_p(\theta,\lambda)$ defined in (3.7). Recall that from Theorem 2, this term is, up to some constant and residual terms, a risk hull. Then denotes by $\theta^\star$ the estimator associated to the filter:

$$\lambda^\star = \arg\min_{\lambda\in\Lambda^\star} U_p(y,\lambda), \tag{4.1}$$

The solution of (4.1) is the penalized blockwise Stein's estimator introduced in (2.6).

**Theorem 4.** *Assume that Assumption A2 holds. Then, there exists $C^\star > 0$ independent of $\epsilon$ such that, for all $\theta \in \ell^2$ and any $0 < \epsilon < 1$:*

$$\mathbb{E}_\theta \|\theta^\star - \theta\|^2 \le (1+\tau_\epsilon) \inf_{\lambda\in\Lambda^\star} R(\theta,\lambda) + C^\star \epsilon^2,$$

*where $\tau_\epsilon \to 0$ as $\epsilon \to 0$ provided $\max_j \mathrm{pen}_j/\sigma_j^2 \to 0$ and $\rho_\epsilon \to 0$ as $\epsilon \to 0$.*

In some sense, this result generalizes Theorem 1. However, we construct a different proof through the risk hull point of view (see Section 6).

Theorem 4 provides in fact an admissible range for the penalty. If we want a sharp oracle inequality, necessarily $\max_j \mathrm{pen}_j/\sigma_j^2 \to 0$ as $\epsilon \to 0$. Hence, the penalty should not be too large. In the same time, we require from Assumption A2 that the penalty contains in a certain sense the variables $(\eta_j)_{j=1\dots J}$. The following lemma provides upper bounds on the term $\mathbb{E}_\theta[\eta_j - \mathrm{pen}_j]_+$ and makes more explicit the behavior of the penalty.

**Lemma 1.** *For all $j \in \{1,\dots,J\}$ and $\delta$ such that $0 < \delta < \epsilon^{-2} b_{K_j-1}^2/2$:*

$$\mathbb{E}[\eta_j - \mathrm{pen}_j]_+ \le \delta^{-1} \exp\left\{ -\delta \mathrm{pen}_j + \delta^2 \Sigma_j^2 + 4\delta^3 \sum_{k\in I_j} \frac{\epsilon^6 b_k^{-6}}{(1 - 2\delta\epsilon^2 b_{K_j-1}^{-2})} \right\},$$

*with*

$$\Sigma_j^2 = \epsilon^4 \sum_{k\in I_j} b_k^{-4}. \tag{4.2}$$

PROOF. Let $j \in \{1, \dots, J\}$ be fixed. First remark that for all $\delta > 0$:

$$\begin{aligned} \mathbb{E}_\theta[\eta_j - \mathrm{pen}_j]_+ &\leq \mathbb{E}_\theta \eta_j \mathbf{1}_{\{\eta_j \geq \mathrm{pen}_j\}} \\ &= \int_{\mathrm{pen}_j}^{+\infty} P(\eta_j \geq t) dt, \\ &= \int_{\mathrm{pen}_j}^{+\infty} P(\exp(\delta\eta_j) \geq e^{\delta t}) dt, \\ &\leq \delta^{-1} e^{-\delta \mathrm{pen}_j} \mathbb{E}_\theta \exp(\delta \eta_j). \end{aligned}$$

Then, provided $0 < \delta < \epsilon^{-2} b_{K_{j-1}}^2 / 2$:

$$\mathbb{E}_\theta \exp(\delta\eta_j) \leq \exp\left\{ \delta^2 \Sigma_j^2 + 4\delta^3 \sum_{k \in I_j} \frac{\epsilon^6 b_k^{-6}}{(1 - 2\epsilon^2 \delta b_k^{-2})_+} \right\}.$$

This conclude the proof of Lemma 1.

□

Let $u$ and $v$ two real sequences. Here and in the sequel, for all $k \in \mathbb{N}$, we write $u_k \lesssim v_k$ if we can find a positive constant $C$ independent of $k$ such that $u_k \leq C v_k$, and $u_k \simeq v_k$ if both $u_k \lesssim v_k$ and $u_k \gtrsim v_k$.

From Assumption A2 and Lemma 1, it is possible to prove that $(\mathrm{pen}_j)_{j=1\dots J}$ should at least fulfill $\mathrm{pen}_j \gtrsim \Sigma_j$ for all $j \in \{1, \dots, J\}$. Since we require in the same time $\max_j \sigma_j^2 / \mathrm{pen}_j \to 0$ as $j \to +\infty$, an admissible penalty in the sense of Theorem 3 should satisfy:

$$\Sigma_j \lesssim \mathrm{pen_j} \lesssim \sigma_j^2, \ \forall\, j \in \{1, \dots, J\}. \tag{4.3}$$

In Section 5, we present a discussion on the choice of the penalty. A specific assumption on $(b_k)_{k \in \mathbb{N}}$ make the situation easier to handle. We present examples where $\mathrm{pen}_j / \Sigma_j \to +\infty$ as $j \to +\infty$ and $\epsilon \to 0$. In particular, we will see that the penalty $\mathrm{pen}_j = \varphi_j \sigma_j^2$ with $\varphi_j = \Delta_j^\gamma$ and $0 < \gamma < 1/2$ is admissible following (4.3).

We are now able to propose an oracle inequality on $\Lambda_{mon}$, the family of monotone filters.

**Corollary 1.** *Assume that Assumption A2 holds and*

$$\max_{j=1\dots J-1} \frac{\sigma_{j+1}^2}{\sigma_j^2} \leq 1 + \eta_\epsilon, \text{ for } 0 < \eta_\epsilon < 1/2. \tag{4.4}$$

*Let $r \geq 0$ be fixed and $N \geq \max\{m : \sum_{k=1}^m b_k^{-2} \leq r^2 \epsilon^{-2} \eta_\epsilon^{-2}\}$. Then , for any $\theta \in \ell^2$ such that $\|\theta\| \leq r$, we have:*

$$\mathbb{E}_\theta \|\theta^\star - \theta\|^2 \leq (1 + \Gamma_\epsilon) \inf_{\lambda \in \Lambda_{mon}} R(\theta, \lambda) + C_3 \epsilon^2,$$

*where $C_3$ denotes a positive constant independent of $\epsilon$, $\Gamma_\epsilon = (2\eta_\epsilon + \tau_\epsilon)/(1-2\eta_\epsilon)$ and $\tau_\epsilon$ is introduced in Theorem 4.*

The proof is a direct consequence of Lemma 1 of (8).

## 5. The choice of the penalty

Following Theorems 2 and 4, the penalized blockwise Stein's estimator is derived from the risk hull method. In the previous section, we have presented conditions on the penalty $(\mathrm{pen}_j)_j$ that lead to sharp oracle inequalites. In this context, the following question naturally arises: what is the smaller possible hull, i.e. may the criterion constructed in (8) be refined?

The Stein hull $V(\theta, \lambda)$ being, up to some constants, the sum of the quadratic risk and a penalty, the previous question may be transfered to the penalty: how can we choose this quantity? This question is quite difficult. The answer cannot be summarized in a single paper. This section only try to shed some light on the properties of the penalty.

Concerning the sequence $(b_k)_{k\in\mathbb{N}}$, we set the following condition:

**Assumption B1**: *There exists $\beta > 0$ such that $(b_k)_{k\in\mathbb{N}} \simeq (k^{-\beta})_{k\in\mathbb{N}}$, i.e. there exist $\bar{b}$ and $\tilde{b}$ independent of $k$ such that $\bar{b}k^{-\beta} \leq b_k \leq \tilde{b}k^{-\beta}$, for all $k \in \mathbb{N}$.*

The eigenvalues are supposed to be polynomially decreasing. The problem is said to be mildly ill-posed.

Assume that B1 holds. Then, for all $j \in \{1 \dots J\}$:

$$\sigma_j^2 = \epsilon^2 \sum_{k\in I_j} b_k^{-2} \simeq \epsilon^2 K_j^{2\beta}(K_{j+1} - K_j),$$

and

$$\Delta_j = \frac{K_{j+1}^{2\beta}}{K_j^{2\beta}(K_{j+1} - K_j)} \simeq (K_{j+1} - K_j)^{-1},$$

provided $K_{j+1}/K_j \to 1$ as $j \to +\infty$. We deduced from Theorem 2 in Section 3 an admissible range (in the oracle sense) for the penalty:

$$\Sigma_j \lesssim \mathrm{pen}_j \lesssim \sigma_j^2, \ \forall j \in \{1, \dots, J\},$$

which, following Assumptions A3 is equivalent to:

$$\epsilon^2 K_j^{2\beta}(K_{j+1} - K_j)^{1/2} \lesssim \mathrm{pen}_j \lesssim \epsilon^2 K_{j+1}^{2\beta}(K_{j+1} - K_j), \ \forall j \in \{1, \dots, J\}.$$

Since, all the penalties in this admissible range may lead to a sharp oracle inequality, one may be interested in finding the smaller possible one, or at least as close as possible to the lower bound of this range.

For all $j \in \{1, \dots, J\}$, the penalty of (8) is:

$$\text{pen}_j^{CT} = \varphi_j \sigma_j^2 = \Delta_j^\gamma \sigma_j^2 \simeq \epsilon^2 K_j^{2\beta} (K_{j+1} - K_j)^{1-\gamma}, \text{ with } 0 < \gamma < 1/2.$$

For $\gamma = 1/2$, Assumption A1 does not hold. The penalized blockwise Stein's estimator is not sharp following (1.6). Remark that this particular choice exactly corresponds to the lower bound of the range (4.3).

The principle of risk hull minimization may lead to more accurate choices. The only restriction on $(\text{pen}_j)_j$ from the risk hull point of view is expressed through Assumption A2:

$$\sum_{j=1}^{J} \mathbb{E}_\theta \left[ \eta_j - \text{pen}_j \right]_+ \leq C_2 \epsilon^2,$$

for some positive constant $C_2$. Since $\mathbb{E}_\theta \left[ \eta_j - 2u \right]_+ \leq \mathbb{E}_\theta \eta_j \mathbf{1}_{\{\eta_j \geq u\}}$ for all positive $u$, we may be interested in the penalty:

$$\overline{\text{pen}}_j = (1+\alpha) \inf \left\{ u : \mathbb{E}_\theta \eta_j \mathbf{1}_{\{\eta_j \geq u\}} \leq \epsilon^2 \right\}, \ \forall j \in \{1, \dots, J\}. \tag{5.1}$$

for some positive $\alpha > 0$.

In the following, we prove that the sequence $(\overline{\text{pen}}_j)_{j=1\dots J}$ defined in (5.1) leads to a sharp oracle inequality. In particular, Assumption A2 holds, i.e. the penalty contains the variability of the problem. For the sake of convenience, we restrict ourselves to one specific type of blocks. All the results presented in the sequel hold for other constructions (see for instance (7)). We leave the proof to the interested reader.

Let $\nu_\epsilon = \lceil \log \epsilon^{-1} \rceil$ and $\kappa_\epsilon = \log^{-1} \nu_\epsilon$, where for all $x \in \mathbb{R}$, $\lceil x \rceil$ denotes the minimal integer strictly greater than $x$. Define the sequence $(T_j)_{j=1\dots J}$ by:

$$\begin{cases} T_1 = \lceil \nu_\epsilon \rceil, \\ T_j = \lceil \nu_\epsilon (1+\kappa_\epsilon)^{j-1} \rceil, \quad j > 1, \end{cases} \tag{5.2}$$

and the bandwidth $J$ as:

$$J = \min\{j : K_j > \bar{N}\}, \text{ with } \bar{N} = \max \left\{ m : \sum_{k=1}^{m} b_k^{-2} \leq \epsilon^{-2} \kappa_\epsilon^{-3} \right\}. \tag{5.3}$$

**Assumption B2**: *The length of the blocks is defined by the sequence $(T_j)_{j=1\dots J}$ where $J$ satisfies (5.3) and the terms $T_j$ are defined in (5.2). For all $j \in \{1, \dots, J\}$, the penalty is $\overline{\text{pen}}_j = (1+\alpha) U_j$, for some $\alpha > 0$.*

In practice, the penalty (5.1) can be computed using Monte-Carlo approximations. We may also use the following lower bound:

**Lemma 2.** *Assume that Assumptions B1 and B2 hold. Then there exists a constant $C$ independent of $j$ and $\epsilon$ such that:*

$$U_j = \inf\left\{u : \mathbb{E}_\theta \eta_j \mathbf{1}_{\{\eta_j \geq u\}} \leq \epsilon^2\right\} \geq \sqrt{2\Sigma_j^2 \log\left(C\epsilon^{-4}\Sigma_j^2\right)}, \tag{5.4}$$

*where for all $j \in \{1 \dots J\}$, $\Sigma_j^2$ is defined in (4.2).*

The proof can be directly derived from the Lemma 1 of (5). Now, if Assumption A3 holds:

$$\overline{\text{pen}}_j \gtrsim K_j^{2\beta}(K_{j+1} - K_j)^{1/2}\sqrt{\log(C\Sigma_j^2)}, \ \forall j \in \{1, \dots, J\}.$$

This penalty reaches, up to a log term, the lower bound of the range (4.3). As we will see the related hull $V(\theta, \lambda)$ is close to the smaller possible one, provided Assumption B1 holds.

Next corollary establishes that the sequence $((1+\alpha)U_j)_{j=1\dots J}$, where the $U_j$ are defined in (5.4), is a relevant choice for the penalty.

**Corollary 2.** *Let $\theta^\star$ the estimator introduced in (2.6). Assume that Assumption B1 and B2 hold. Then,*

$$\mathbb{E}_\theta\|\theta^\star - \theta\|^2 \leq (1 + \gamma_\epsilon) \inf_{\lambda \in \Lambda^\star} R(\theta, \lambda) + \frac{C_4}{\alpha}\epsilon^2,$$

*where $C_4$ denotes a positive constant independent of $\epsilon$ and $\gamma_\epsilon = o(1)$ as $\epsilon \to 0$.*

PROOF. From Theorems 2 and 4, we only have to prove that Assumption A2 holds since $\max_j \overline{\text{pen}}_j/\sigma_j^2$ converges to 0 as $\epsilon \to 0$. For all $j \in \{1, \dots, J\}$, using Lemma 1 in Section 4:

$$\mathbb{E}\left[\eta_j - \overline{\text{pen}}_j\right]_+ \leq \frac{1}{\delta}\exp\left\{-\delta\overline{\text{pen}}_j + \delta^2\Sigma_j^2 + 4\delta^3 \sum_{k \in I_j} \frac{\epsilon^6 b_k^{-6}}{(1 - 2\delta\epsilon^2 b_{K_{j-1}}^{-2})}\right\}, \tag{5.5}$$

for all $0 < \delta < \epsilon^{-2} b_{K_j - 1}^2/2$. Setting,

$$\delta = \sqrt{\frac{\log(C\epsilon^{-4}\Sigma_j^2)}{2\Sigma_j^2}},$$

and using Lemma 2, we obtain:

$$\begin{aligned}
&\mathbb{E}[\eta_j - \overline{\text{pen}}_j]_+ \\
&\leq \sqrt{\frac{2\Sigma_j^2}{\log(C\epsilon^{-4}\Sigma_j^2)}} \exp\left\{\frac{1}{2}\log(C\epsilon^{-4}\Sigma_j^2)\right\} \times \exp\left\{-(1+\alpha)\log(C\epsilon^{-4}\Sigma_j^2)\right\}, \\
&\leq C\epsilon^2\sqrt{\frac{1}{\log(C\epsilon^{-4}\Sigma_j^2)}} \exp\left\{-\alpha\log(C\epsilon^{-4}\Sigma_j^2)\right\}.
\end{aligned}$$

Indeed, provided (5.2), (5.3) and Assumption B1 hold, $\delta b_{K_j-1}^{-2}$ and the last term in the right hand side of the exponential in (5.5) converge to 0 as $j \to +\infty$. Hence, we eventually obtain:

$$\begin{aligned}
\sum_{j=1}^{J} \mathbb{E}[\eta_j - \overline{\text{pen}}_j]_+ &\leq C\epsilon^2 \sum_{j=1}^{J} \frac{1}{\log^{1/2}(CT_j)} \exp\{-\alpha \log(CT_j)\}, \\
&\leq C\epsilon^2 \sum_{j=1}^{J} j^{-1/2} \exp\left\{-\alpha \log(C\nu_\epsilon(1+\kappa_\epsilon)^j)\right\}, \\
&\leq C\epsilon^2 \sum_{j=1}^{+\infty} j^{-1/2} \exp\{-\alpha D j\} < \frac{C\epsilon^2}{\alpha},
\end{aligned}$$

where $D$ and $C$ denote two positive constants independent of $\epsilon$. This concludes the proof of Corollary 2.

□

From Corollary 2, the penalty $(\overline{\text{pen}}_j)_{j=1,\dots,J}$ leads to a sharp oracle inequality. Using the same bounds as in the proof, it seems hopeless to obtain a similar result with the penalty $(\Sigma_j)_{j=1\dots J}$. Remark that with our assumptions, this choice exactly corresponds to the penalty of (8) with $\gamma = 1/2$. From (5.5), the Assumption A2 will not be satisfied. However, the bound proposed by the Lemma 1 may be perhaps improved or this Assumption A2 relaxed.

In a certain sense, the risk hull point of view quantifies the effects of the penalty. We have presented an admissible range for this quantity and proposed a choice close to the lower bound of this range. From now, the question is: 'What is a good penalty?'. A small penalty leads to a sharp hull (see Theorem 2) but the constant $C^\star$ in Theorem 3 may be large. On the opposite hand, large penalties perfectly contain the variability of the problem but the hull is less precise. Hence, it is not clear that such a problem may be easily solved. This was not the goal of the present paper.

In order to conclude, it seems necessary to discuss about the role played by the constant $\alpha$ in the penalty $(\overline{\text{pen}}_j)_{j=1\dots J}$. Assumption A2 does not hold for $\alpha = 0$. On the other hand, the proof of Corollary 2 indicates that large values for $\alpha$ will not lead to an accurate recovering. The choice of $\alpha$ has already been discussed and illustrated via some numerical simulations in a sligthy different setting: see (5) or (17) for more details. Remark that we do not require $\alpha$ to be greater than 1 in this paper. This is a small difference with the constraints expressed in a regularization parameter choice scheme. This can be explained by the blockwise structure of the variables $(\eta_j)_{j=1\dots J}$.

## 6. Proofs and technical lemmas

### 6.1. Ordered processes

Ordered processes were introduced in (16). The more recent paper (4) studies in detail these processes and provides very interesting tools. These stochastic objects may play an important role in adaptive estimation: see in particular (13) or (17) for more details.

The aim of this section is not to provide an exhaustive presentation of this theory but rather to introduce some definitions and useful properties.

**Definition 1.** *Let $\zeta(t)$, $t \geq 0$ a separable random process with $\mathbb{E}\zeta(t) = 0$ and finite variance $\Sigma^2(t)$. It is called ordered if for all $t_2 \geq t_1 \geq 0$,*

$$\Sigma^2(t_2) \geq \Sigma^2(t_1), \text{ and } \mathbb{E}[\zeta(t_2) - \zeta(t_1)]^2 \leq \Sigma^2(t_2) - \Sigma^2(t_1). \tag{6.1}$$

Let $\zeta$ a standard Gaussian random variable. The process $t \mapsto \zeta t$ is the most simple example of ordered process. Wiener processes are also covered by Definition 1. The family of ordered processes is in fact quite large.

**Assumption C1**. *There exists $\kappa > 0$ such that:*

$$\varphi(\kappa) = \sup_{t_1,t_2} \mathbb{E} \exp\left\{ \kappa \frac{\zeta(t_1) - \zeta(t_2)}{\sqrt{\mathbb{E}[\zeta(t_1) - \zeta(t_2)]^2}} \right\} < +\infty. \tag{6.2}$$

This assumption is not very restrictive. Several processes encountered in linear estimation satisfy this hypothesis.

The proof of the following result can be found in (4).

**Lemma 3.** *Let $\zeta(t)$, $t \geq 0$ an ordered process satisfying $\zeta(0) = 0$ and Assumption C1. There exists a constant $C = C(\kappa)$ such that for all $\gamma > 0$:*

$$\mathbb{E} \sup_{t\geq 0} \left[\zeta(t) - \gamma\Sigma^2(t)\right]_+ \leq \frac{C}{\gamma}.$$

This lemma is rather important in the theory of ordered processes and leads to several interesting results. In particular, the following corollary will be often used in the proofs.

**Corollary 3.** *Let $\zeta(t)$, $t \geq 0$ an ordered process satisfying $\zeta(0) = 0$ and Assumption C1. Consider $\hat{t}$ measurable with respect to $\zeta$. Then, there exists $C = C(\kappa) > 0$ such that:*

$$\mathbb{E}\zeta(\hat{t}) \leq C\sqrt{\mathbb{E}\Sigma^2(\hat{t})}.$$

PROOF. Let $\gamma > 0$ be fixed. Using Lemma 1,

$$\begin{aligned}
\mathbb{E}\zeta(\hat{t}) &= \mathbb{E}\zeta(\hat{t}) - \gamma\mathbb{E}\Sigma^2(\hat{t}) + \gamma\mathbb{E}\Sigma^2(\hat{t}), \\
&\leq \mathbb{E}\sup_{t\geq 0}\left[\zeta(t) - \gamma\Sigma^2(t)\right]_+ + \gamma\mathbb{E}\Sigma^2(\hat{t}), \\
&\leq \frac{C}{\gamma} + \gamma\mathbb{E}\Sigma^2(\hat{t})
\end{aligned}$$

Choose $\gamma = \sqrt{\mathbb{E}\Sigma^2(\hat{t})}$ in order to conclude the proof.

□

### 6.2. Proofs of Theorems 2-4

**Proof of Theorem 2**. First, remark that:

$$\begin{aligned}
&\mathbb{E}_\theta \sup_{\lambda\in\Lambda^\star} \left\{ \|\hat{\theta}_\lambda - \theta\|^2 - V(\theta,\lambda) \right\} \\
&= \mathbb{E}_\theta \sup_{\lambda\in\Lambda^\star} \left\{ \sum_{k=1}^{+\infty} (1-\lambda_k)^2\theta_k^2 + \epsilon^2 \sum_{k=1}^{+\infty} \lambda_k^2 b_k^{-2}\xi_k^2 - 2\epsilon \sum_{k=1}^{+\infty} \lambda_k(1-\lambda_k)\theta_k b_k^{-1}\xi_k \right. \\
&\qquad \left. - V(\theta,\lambda) \right\}, \\
&= \mathbb{E}_\theta \sup_{\lambda\in\Lambda^\star} \left\{ \sum_{j=1}^{J} \left[ (1-\lambda_j)^2\|\theta\|_{(j)}^2 + \lambda_j^2 \sum_{k\in I_j} \epsilon^2 b_k^{-2}\xi_k^2 - 2\lambda_j(1-\lambda_j)X_j \right] \right. \\
&\qquad \left. + \sum_{k>N} \theta_k^2 - V(\theta,\lambda) \right\}, \\
&= \mathbb{E}_\theta \sum_{j=1}^{J} \left[ (1-\hat{\lambda}_j)^2\|\theta\|_{(j)}^2 + \hat{\lambda}_j^2 \sum_{k\in I_j} \epsilon^2 b_k^{-2}\xi_k^2 + 2\hat{\lambda}_j(\hat{\lambda}_j-1)X_j \right] \\
&\qquad + \sum_{k>N} \theta_k^2 - \mathbb{E}_\theta V(\theta,\hat{\lambda}),
\end{aligned}$$

with

$$\hat{\lambda} = \arg\sup_{\lambda\in\Lambda^\star} \left\{ \|\hat{\theta}_\lambda - \theta\|^2 - V(\theta,\lambda) \right\},$$

and

$$X_j = \epsilon \sum_{k\in I_j} \theta_k b_k^{-1}\xi_k, \ \forall j \in \{1,\dots,J\}. \tag{6.3}$$

Let $j \in \{1,\dots,J\}$ be fixed. Use the decomposition:

$$\begin{aligned}
\mathbb{E}_\theta 2\hat{\lambda}_j(\hat{\lambda}_j-1)X_j &= \mathbb{E}_\theta \hat{\lambda}_j^2 X_j + \mathbb{E}_\theta(\hat{\lambda}_j^2 - 2\hat{\lambda}_j)X_j, \\
&= \mathbb{E}_\theta \hat{\lambda}_j^2 X_j + \mathbb{E}_\theta(1-\hat{\lambda}_j)^2 X_j = A_j^1 + A_j^2,
\end{aligned} \tag{6.4}$$

since $\mathbb{E}_\theta X_j = 0$. First consider $A_j^1$. Let $\lambda_j^0$ denotes the blockwise constant oracle on the block $j$. Using Corollary 3 in Section 6.1:

$$A_j^1 = \mathbb{E}_\theta \hat{\lambda}_j^2 X_j = \mathbb{E}_\theta \left[ \hat{\lambda}_j^2 - (\lambda_j^0)^2 \right] X_j \le C \sqrt{\mathbb{E}_\theta \left[ \hat{\lambda}_j^2 - (\lambda_j^0)^2 \right]^2 \sum_{k\in I_j} \epsilon^2 b_k^{-2}\theta_k^2}, \tag{6.5}$$

where $C > 0$ denotes a positive constant. Indeed, both processes $\zeta : t \mapsto (t^2 - (\lambda_j^0)^2)X_j,\ t \in [(\lambda_j^0)^2, 1]$ and $\bar{\zeta} : t \mapsto (t^{-2} - (\lambda_j^0)^2)X_j,\ t \in [(\lambda_j^0)^{-1}; +\infty[$ are ordered and satisfy Assumption C1. For all $\gamma > 0$, using:

$$\left[\hat{\lambda}_j^2 - (\lambda_j^0)^2\right]^2 \leq 4\left[(1-\hat{\lambda}_j)^2 + (1-\lambda_j^0)^2\right](\hat{\lambda}_j^2 + (\lambda_j^0)^2), \tag{6.6}$$

and the Cauchy-Schwartz and Young inequalities:

$$\begin{aligned} A_j^1 &\leq C\sqrt{\mathbb{E}_\theta\left[(1-\hat{\lambda}_j)^2 + (1-\lambda_j^0)^2\right](\hat{\lambda}_j^2 + (\lambda_j^0)^2)\max_{k\in I_j}\epsilon^2 b_k^{-2}\|\theta\|_{(j)}^2}, \\ &\leq C\mathbb{E}_\theta\left[\gamma(1-\hat{\lambda}_j)^2\|\theta\|_{(j)}^2 + \gamma^{-1}\Delta_j\hat{\lambda}_j^2\sigma_j^2\right] + C\gamma(1-\lambda_j^0)^2\|\theta\|_{(j)}^2 \\ &\quad + C\gamma^{-1}\Delta_j(\lambda_j^0)^2\sigma_j^2 + C\sqrt{\mathbb{E}_\theta(1-\hat{\lambda}_j)^2\hat{\lambda}_j^2\max_{k\in I_j}\epsilon^2 b_k^{-2}\|\theta\|_{(j)}^2}, \end{aligned} \tag{6.7}$$

for some positive constant $C$. The bound of the last term in the r.h.s. of (6.7) requires to be careful. In a first time, suppose that:

$$\|\theta\|_{(j)}^2 \leq \sigma_j^2. \tag{6.8}$$

In such a situation, for all $\gamma > 0$:

$$\begin{aligned} \sqrt{\mathbb{E}_\theta(1-\hat{\lambda}_j)^2\hat{\lambda}_j^2\max_{k\in I_j}\epsilon^2 b_k^{-2}\|\theta\|_{(j)}^2} &\leq \sqrt{\|\theta\|_{(j)}^2\mathbb{E}_\theta\hat{\lambda}_j^2\max_{k\in I_j}\epsilon^2 b_k^{-2}}, \\ &\leq \gamma\|\theta\|_{(j)}^2 + \gamma^{-1}\Delta_j\mathbb{E}_\theta\hat{\lambda}_j^2\sigma_j^2. \end{aligned}$$

If (6.8) holds, then:

$$\|\theta\|_{(j)}^2 = \frac{\sigma_j^2\|\theta\|_{(j)}^2}{\sigma_j^2 + \|\theta\|_{(j)}^2}\left(1 + \frac{\|\theta\|_{(j)}^2}{\sigma_j^2}\right) \leq 2\left\{(1-\lambda_j^0)^2\|\theta\|_{(j)}^2 + (\lambda_j^0)^2\sigma_j^2\right\}, \tag{6.9}$$

where $\lambda^0$ is the oracle defined in (2.2). Indeed,

$$\lambda_j^0 = \frac{\|\theta\|_{(j)}^2}{\sigma_j^2 + \|\theta\|_{(j)}^2},\ \forall j \in \{1, \dots, J\}.$$

Now, suppose:

$$\|\theta\|_{(j)}^2 > \sigma_j^2. \tag{6.10}$$

Then, for all $\gamma > 0$:

$$\begin{aligned} \sqrt{\mathbb{E}_\theta(1-\hat{\lambda}_j)^2\hat{\lambda}_j^2\max_{k\in I_j}\epsilon^2 b_k^{-2}\|\theta\|_{(j)}^2} &\leq \sqrt{\max_{k\in I_j}\epsilon^2 b_k^{-2}\mathbb{E}_\theta(1-\hat{\lambda}_j)^2\|\theta\|_{(j)}^2}, \\ &\leq \gamma\mathbb{E}_\theta(1-\hat{\lambda}_j)^2\|\theta\|_{(j)}^2 + \gamma^{-1}\Delta_j\sigma_j^2. \end{aligned}$$

Using (6.10):

$$\sigma_j^2 = \frac{\sigma_j^2\|\theta\|_{(j)}^2}{\sigma_j^2 + \|\theta\|_{(j)}^2}\left(1 + \frac{\sigma_j^2}{\|\theta\|_{(j)}^2}\right) \leq 2\left\{(1-\lambda_j^0)^2\|\theta\|_{(j)}^2 + (\lambda_j^0)^2\sigma_j^2\right\}.$$

Setting $\gamma = \sqrt{\Delta_j}$, we eventually obtain:

$$A_j^1 \leq C\sqrt{\Delta_j}\mathbb{E}_\theta\left[(1-\hat{\lambda}_j)^2\|\theta\|^2_{(j)} + \hat{\lambda}_j^2\sigma_j^2\right] + C\sqrt{\Delta_j}\left[(1-\lambda_j^0)^2\|\theta\|^2_{(j)} + (\lambda_j^0)^2\sigma_j^2\right], \tag{6.11}$$

for some constant $C>0$ independent of $\epsilon$. The same bound occurs for the term $A_j^2$ in (6.4). Hence, there exists $B>0$ independent of $\epsilon$ such that:

$$\begin{aligned}
&\mathbb{E}_\theta \sup_{\lambda\in\Lambda^\star}\left\{\|\hat{\theta}_\lambda - \theta\|^2 - V(\theta,\lambda)\right\}\\
&\leq \mathbb{E}_\theta \sum_{j=1}^J \left[(1+B\rho_\epsilon)(1-\hat{\lambda}_j)^2\|\theta\|^2_{(j)} + \hat{\lambda}_j^2\sum_{k\in I_j}\epsilon^2 b_k^{-2}\xi_k^2 + B\rho_\epsilon\hat{\lambda}_j^2\sigma_j^2\right]\\
&\qquad + \sum_{k>N}\theta_k^2 + B\rho_\epsilon R(\theta,\lambda^0) - \mathbb{E}_\theta V(\theta,\hat{\lambda}),\\
&\leq \mathbb{E}_\theta \sup_{\lambda\in\Lambda^\star}\left\{\sum_{j=1}^J\left[(1+B\rho_\epsilon)(1-\lambda_j)^2\|\theta\|^2_{(j)} + \lambda_j^2\sum_{k\in I_j}\epsilon^2 b_k^{-2}\xi_k^2 + B\rho_\epsilon\lambda_j^2\sigma_j^2\right]\right.\\
&\qquad \left. + \sum_{k>N}\theta_k^2 + B\rho_\epsilon R(\theta,\lambda^0) - V(\theta,\lambda)\right\},
\end{aligned}$$

where $\rho_\epsilon$ is defined in (3.4). Now, using (3.3) and (3.6),

$$\begin{aligned}
&\mathbb{E}_\theta \sup_{\lambda\in\Lambda^\star}\left\{\|\hat{\theta}_\lambda - \theta\|^2 - V(\theta,\lambda)\right\}\\
&\leq \mathbb{E}_\theta \sup_{\lambda\in\Lambda^\star}\left\{\sum_{j=1}^J\left[\lambda_j^2\eta_j - 2\lambda_j \mathrm{pen}_j\right] - C_2\epsilon^2\right\},\\
&= \sum_{j=1}^J \mathbb{E}_\theta \sup_{\lambda_j\in[0,1]}\left[\lambda_j^2\eta_j - 2\lambda_j\mathrm{pen}_j\right] - C_2\epsilon^2.
\end{aligned}$$

Let $j\in\{1\dots J\}$ be fixed. We are looking for $\lambda_j\in[0,1]$ that maximizes the quantity $\lambda_j^2\eta_j - 2\lambda_j\mathrm{pen}_j$. If $\eta_j<0$, the function $\lambda\mapsto\lambda^2\eta_j - 2\lambda\mathrm{pen}_j$ is concave and the maximum on $[0,1]$ is attained for $\lambda=0$. Now, if $\eta_j>0$, the function $\lambda\mapsto\lambda^2\eta_j - 2\lambda\mathrm{pen}_j$ is convex and the maximum on $[0,1]$ is attained in 0 or in 1. Therefore:

$$\sup_{\lambda_j\in[0,1]}\left\{\lambda_j^2\eta_j - 2\lambda_j\mathrm{pen}_j\right\} = \left[\eta_j - 2\mathrm{pen}_j\right]_+, \tag{6.12}$$

Using Assumption A2, we eventually obtain:

$$\begin{aligned}
\mathbb{E}_\theta \sup_{\lambda\in\Lambda^\star}\left\{\|\hat{\theta}_\lambda - \theta\|^2 - V(\theta,\lambda)\right\} &\leq \sum_{j=1}^J \mathbb{E}_\theta\left[\eta_j - 2\mathrm{pen}_j\right]_+ - C_2\epsilon^2,\\
&\leq \sum_{j=1}^J \mathbb{E}_\theta\left[\eta_j - \mathrm{pen}_j\right]_+ - C_2\epsilon^2 \leq 0.
\end{aligned}$$

This concludes the proof of Theorem 2.

$\square$

**Proof of Theorem 3**. Using the same algebra of the proof of Theorem 2, we obtain:

$$\mathbb{E}_\theta \sup_{\lambda\in\Lambda^\star}\left\{\|\hat\theta_\lambda-\theta\|^2-W(\theta,\lambda)\right\}\leq \sum_{j=1}^J \mathbb{E}_\theta \sup_{\lambda_j\in[0,1]}\{\lambda_j^2\eta_j-\lambda_j^2\mathrm{pen}_j\}-C_2\epsilon^2.$$

For all $j\in\{1\dots J\}$, the only difference with $V(\theta,\lambda)$ is contained in the bound of:

$$\sup_{\lambda_j\in[0,1]}\left\{\lambda_j^2\eta_j-\lambda_j^2\mathrm{pen}_j\right\}\leq[\eta_j-\mathrm{pen}_j]_+.$$

Using Assumption A2, we obtain Theorem 3.

$\square$

**Proof of Theorem 4**. In the situation where Assumption A2 holds, Theorem 3 provides that:

$$\mathbb{E}_\theta\|\theta^\star-\theta\|^2\leq W(\theta,\lambda^\star)=(1+B\rho_\epsilon)\bar R_p(\theta,\lambda^\star)+B\rho_\epsilon R(\theta,\lambda^0)+C_2\epsilon^2, \quad (6.13)$$

where:

$$\bar R_p(\theta,\lambda^\star)=\sum_{j=1}^J\left[(1-\lambda_j^\star)^2\|\theta\|_{(j)}^2+(\lambda_j^\star)^2\sigma_j^2+(\lambda_j^\star)^2\mathrm{pen}_j\right]+\sum_{k>N}\theta_k^2,$$

and $B$ denotes a positive constant independent of $\epsilon$. Moreover, from (4.1),

$$U_p(y,\lambda^\star)\leq U_p(y,\lambda),\ \forall\lambda\in\Lambda^\star.$$

The proof of Theorem 4 is mainly based on this two equalities. First remark that:

$$\begin{aligned}
&U_p(y,\lambda^\star)-\bar R_p(\theta,\lambda^\star)\\
&\quad=\sum_{j=1}^J\Big[\{(\lambda_j^\star)^2-2\lambda_j^\star\}(\|\tilde y\|_{(j)}^2-\sigma_j^2)+(\lambda_j^\star)^2\sigma_j^2+2\lambda_j^\star\mathrm{pen}_j-(1-\lambda_j^\star)^2\|\theta\|_{(j)}^2\\
&\qquad\qquad-(\lambda_j^\star)^2\sigma_j^2-(\lambda_j^\star)^2\mathrm{pen}_j\Big]-\sum_{k>N}\theta_k^2,\\
&\quad=\sum_{j=1}^J\Big[\{(\lambda_j^\star)^2-2\lambda_j^\star\}(\|\tilde y\|_{(j)}^2-\sigma_j^2)-(1-\lambda_j^\star)^2\|\theta\|_{(j)}^2+\{2\lambda_j^\star\\
&\qquad\qquad-(\lambda_j^\star)^2\}\mathrm{pen}_j\Big]-\sum_{k>N}\theta_k^2.
\end{aligned}$$

Hence,

$$U_p(y,\lambda^\star) - \bar{R}_p(\theta,\lambda^\star)$$

$$= \sum_{j=1}^{J} \Bigg[ \{(\lambda_j^\star)^2 - 2\lambda_j^\star\} \sum_{k\in I_j} (\theta_k^2 + \epsilon^2 b_k^{-2}(\xi_k^2 - 1) + 2\epsilon b_k^{-1}\xi_k\theta_k) - (1-\lambda_j^\star)^2 \|\theta\|_{(j)}^2$$

$$+\{2\lambda_j^\star - (\lambda_j^\star)^2\}\mathrm{pen}_j \Big] - \sum_{k>N} \theta_k^2,$$

$$= \sum_{j=1}^{J} \{(\lambda_j^\star)^2 - 2\lambda_j^\star\}(\eta_j + 2X_j - \mathrm{pen}_j) - \|\theta\|^2,$$

where $\eta_j$ and $X_j$ are respectively defined in (3.3) and (6.3). Hence, from (4.1),

$$\bar{R}_p(\theta,\lambda^\star) = U_p(y,\lambda^\star) + \|\theta\|^2 + \sum_{j=1}^{J} \{2\lambda_j^\star - (\lambda_j^\star)^2\}(\eta_j + 2X_j - \mathrm{pen}_j),$$

$$\leq U_p(y,\lambda^p) + \|\theta\|^2 + \sum_{j=1}^{J} \{2\lambda_j^\star - (\lambda_j^\star)^2\}(\eta_j + 2X_j - \mathrm{pen}_j), \quad (6.14)$$

where

$$\lambda^p = \arg\inf_{\lambda\in\Lambda^\star} R_p(\theta,\lambda).$$

and $R_p(\theta,\lambda)$ is defined in (3.7). Then, with simple algebra:

$$\mathbb{E}_\theta U_p(y,\lambda^p) = \mathbb{E}_\theta \sum_{j=1}^{J} \left[ \{(\lambda_j^p)^2 - 2\lambda_j^p\}(\|\tilde{y}\|_{(j)}^2 - \sigma_j^2) + (\lambda_j^p)^2\sigma_j^2 + 2\lambda_j^p \mathrm{pen}_j \right],$$

$$= R_p(\theta,\lambda^p) - \|\theta\|^2.$$

This leads to,

$$\mathbb{E}_\theta \bar{R}_p(\theta,\lambda^\star) \leq R_p(\theta,\lambda^p) + \mathbb{E}_\theta \sum_{j=1}^{J} \{2\lambda_j - (\lambda_j^\star)^2\}(\eta_j + 2X_j - \mathrm{pen}_j). \quad (6.15)$$

We are now interested in the behavior of the right hand side of (6.15). First, using (6.4)-(6.11) in the proof of Theorem 2:

$$\mathbb{E}_\theta\{2\lambda_j - (\lambda_j^\star)^2\}X_j$$

$$\leq C\rho_\epsilon \left\{ (1-\lambda_j^p)^2\|\theta\|_{(j)}^2 + (\lambda_j^p)^2\sigma_j^2 \right\} + \bar{C}\rho_\epsilon \mathbb{E}_\theta \left\{ (1-\lambda_j^\star)^2\|\theta\|_{(j)}^2 + (\lambda_j^\star)^2\sigma_j^2 \right\},$$

for all $j \in \{1,\dots,J\}$. Here, $C$ and $\bar{C}$ denote positive constants independent of $\epsilon$. In particular, it is always possible to obtain $\bar{C}$ verifying $\bar{C}\rho_\epsilon < 1$ (see the proof of Theorem 2 for more details). Hence:

$$\mathbb{E}_\theta \bar{R}_p(\theta,\lambda^\star)$$

$$\leq (1 + C\rho_\epsilon)R_p(\theta,\lambda^p) + \bar{C}\rho_\epsilon \mathbb{E}_\theta \bar{R}_p(\theta,\lambda^\star) + \mathbb{E}_\theta \sum_{j=1}^{J} \{2\lambda_j^\star - (\lambda_j^\star)^2\}(\eta_j - \mathrm{pen}_j).$$

Then, from Assumption A2 and (6.12):

$$\mathbb{E}_\theta \sum_{j=1}^{J} \{2\lambda_j^\star - (\lambda_j^\star)^2\}(\eta_j - \mathrm{pen}_j) = \mathbb{E}_\theta \sum_{j=1}^{J} [\eta_j - \mathrm{pen}_j]_+ \leq C_2 \epsilon^2.$$

This leads to:

$$\begin{aligned} \mathbb{E}_\theta \bar{R}_p(\theta, \lambda^\star) &\leq (1 + C\rho_\epsilon) R_p(\theta, \lambda^p) + \bar{C}\rho_\epsilon \mathbb{E}_\theta \bar{R}_p(\theta, \lambda^\star) + C_2\epsilon^2, \\ \Rightarrow (1 - \bar{C}\rho_\epsilon)\mathbb{E}_\theta \bar{R}_p(\theta, \lambda^\star) &\leq (1 + C\rho_\epsilon) R_p(\theta, \lambda^p) + C_2\epsilon^2, \\ \Rightarrow \mathbb{E}_\theta \bar{R}_p(\theta, \lambda^\star) &\leq \frac{(1 + C\rho_\epsilon)}{(1 - \bar{C}\rho_\epsilon)} R_p(\theta, \lambda^p) + C\epsilon^2. \end{aligned} \tag{6.16}$$

Using (6.13) and (6.16):

$$\begin{aligned} \mathbb{E}_\theta \|\theta^\star - \theta\|^2 &\leq (1 + B\rho_\epsilon)\mathbb{E}_\theta \bar{R}_p(\theta, \lambda^\star) + C_2\epsilon^2 + B\rho_\epsilon R(\theta, \lambda^0), \\ &\leq (1 + \mu_\epsilon) R_p(\theta, \lambda^p) + C\epsilon^2 + B\rho_\epsilon R(\theta, \lambda^0), \end{aligned}$$

where $\mu_\epsilon = \mu_\epsilon(\rho_\epsilon)$ is such that $\mu_\epsilon \to 0$ as $\rho_\epsilon \to 0$ and $C$ is a positive constant independent of $\epsilon$. In order to conclude the proof, we just have to compare $R(\theta, \lambda^0)$ to $R_p(\theta, \lambda_p)$. For all $j \in \{1, \ldots, J\}$, introduce:

$$R_p^j(\theta, \lambda) = (1-\lambda_j)^2 \|\theta\|_{(j)}^2 + \lambda_j^2 \sigma_j^2 + 2\lambda_j \mathrm{pen}_j, \text{ and } R^j(\theta, \lambda) = (1-\lambda_j)^2 \|\theta\|_{(j)}^2 + \lambda_j^2 \sigma_j^2.$$

Then:

$$\begin{aligned} R_p^j(\theta, \lambda^p) &\leq \frac{\sigma_j^4 \|\theta\|_{(j)}^2}{(\sigma_j^2 + \|\theta\|_{(j)}^2)^2} + \frac{\sigma_j^2 \|\theta\|_{(j)}^4}{(\sigma_j^2 + \|\theta\|_{(j)}^2)^2} + 2\frac{\mathrm{pen_j}}{\sigma_j^2} \frac{\sigma_j^2 \|\theta\|_{(j)}^2}{\sigma_j^2 + \|\theta\|_{(j)}^2}, \\ &= \left(1 + 2\frac{\mathrm{pen_j}}{\sigma_j^2}\right) R^j(\theta, \lambda^0), \end{aligned}$$

since $R_p^j(\theta, \lambda^p) \leq R_p^j(\theta, \lambda^0)$ from the definition of $\lambda^p$. This concludes the proof of Theorem 3.

□

## References

[1] Barron, A., Birgé, L. and Massart, P. (1999). Risk bounds for model selection via penalization. *Probability Theory and Related Fields.* **113** 301–413.

[2] Bissantz, N., Hohage, T., Munk, A. and Ryumgaart, F. (2007). Convergence rates of general regularization methods for statistical inverse problems and applications. *SIAM J. Numerical Analysis.* **45** 2610–2636.

[3] Candès, E. (2006). Modern statistical estimation via oracle inequalities. *Acta numerica.* **15** 257–325.

[4] Cao, Y. and Golubev, Y. (2006). On adaptive regression by splines. *Mathematical Methods of Statistics.* **15** 398–414.

[5] Cavalier, L. and Golubev, Y. (2006). Risk hull method and regularization by projections of ill-posed inverse problems. *Ann. Statist.*, **34** 1653–1677.

[6] Cavalier, L., Golubev, Y., Picard, D. and Tsybakov, A.B. (2002). Oracle inequalities for inverse problems. *Annals of Statistics.* **30** 843–874.

[7] Cavalier, L. and Tsybakov, A.B. (2001). Penalized blockwise Stein's method, monotone oracles and sharp adaptative estimation. *Mathematical Methods of Statistics.* **3** 247–282.

[8] Cavalier, L. and Tsybakov, A.B. (2002). Sharp adaptation for inverse problems with random noise. *Probability Theory and Related Fields.* **123** 323–354.

[9] Donoho, D.L. (1995). Nonlinear solutions of linear inverse problems by wavelet-vaguelette decomposition. *Journal of Applied and Computationnal Harmonic Analysis.* **2** 101–126.

[10] Engl, H.W., Hank, M. and Neubauer, A. (1996). Regularization of Inverse Problems. *Kluwer Academic Publishers Group, Dordrecht.*

[11] Ermakov, M.S. (1989). Minimax estimation of the solution of an ill-posed convolution type problem. *Problems of Information Transmission.* **25** 191–200.

[12] Fan, J. (1991). On the optimal rates of convergence for nonparametric deconvolutions problems. *Annals of Statistics.* **19** 1257–1272.

[13] Golubev, Y. (2007). On oracle inequalities related to high dimensional linear models. *IMA Proceedings. 'Topics in stochastic analysis and nonparametric estimation' Chow P.L., Mordukhovich B. and Yin, G. (eds.) Springer-Verlag* 105–122.

[14] Hida, T. (1980). Brownian motion. *Springer-Verlag, New-York.*

[15] Johnstone, I.M. and Silverman, B.W. (1990). Speed of estimation in positron emission tomography and related inverse problems. *Annals of Statistics* **18** 251–280.

[16] Kneip, A. (1994). Ordered linear smoother. *Annals of Statistics.* **22** 835–866.

[17] Marteau, C. (2007). Risk hull method for general families of estimators. *In preparation.*